\documentclass[a4paper, 11pt]{amsart}
\usepackage{amsmath, amssymb, amsthm, amsfonts,enumitem,tikz,relsize}
\usetikzlibrary{calc,positioning}
\usepackage[normalem]{ulem}
\usepackage[colorlinks = true, citecolor = blue, hypertexnames=false]{hyperref}
\usepackage{xcolor}
\hypersetup{
    colorlinks,
    linkcolor={red!70!black},
    citecolor={blue!70!black},
    urlcolor={blue!80!black}
}
\usepackage{framed,float,caption}
\usepackage[OT1]{fontenc}
\usepackage{comment}

\theoremstyle{plain}
\newtheorem{theorem}{Theorem}[section]
\newtheorem{corollary}[theorem]{Corollary}
\newtheorem{lemma}[theorem]{Lemma}
\newtheorem{proposition}{Proposition}
\newtheorem{fact}[theorem]{Fact}

\theoremstyle{definition}
\newtheorem*{definition*}{Definition}
\newtheorem*{example*}{Example}
\newtheorem*{remark}{Remark}

\newtheorem{problem}{Problem}
\newtheorem*{problem*}{Problem}
\newtheorem*{mainproblem}{\hypertarget{h:problem}{Main Problem}}
\newtheorem*{maintheorem}{\hypertarget{h:theorem}{Theorem}}
\newtheorem*{maincorollary}{\hypertarget{h:corollary}{Corollary}}
\newtheorem*{conjecture*}{Conjecture}
\newtheorem*{goal*}{Goal}

\newtheorem*{setup}{Setup}

\theoremstyle{remark}
\newtheorem{claim}{Claim}

\newenvironment{proofclaim}[1][Proof of Claim]{\begin{proof}[#1]}{\end{proof}}

\DeclareMathOperator{\rk}{rk}
\DeclareMathOperator{\gtd}{gtd}
\DeclareMathOperator{\charac}{char}

\DeclareMathOperator{\Sym}{Sym} 
\DeclareMathOperator{\Alt}{Alt}
\DeclareMathOperator{\Inn}{Inn}

\DeclareMathOperator{\GL}{GL}
\DeclareMathOperator{\PGL}{PGL} 
\DeclareMathOperator{\PSL}{PSL}
\DeclareMathOperator{\pSL}{(P)SL}
\DeclareMathOperator{\AGL}{AGL} 
\DeclareMathOperator{\SL}{SL}

\DeclareMathOperator{\pr}{pr} 
\DeclareMathOperator{\mr}{m} 

\newcommand{\cibo}{\textbf{CiBo}}

\begin{document}
\title[Degree of generic transitivity]{Bounding the degree of generic sharp transitivity}
\author{Tuna Alt\i{}nel}
\address{Universit\'e de Lyon \\
Universit\'e Claude Bernard Lyon 1\\
CNRS UMR 5208\\
Institut Camille Jordan\\
43 blvd du 11 novembre 1918\\
F-69622 Villeurbanne cedex\\
France}
\email{altinel@math.univ-lyon1.fr}
\author{Joshua Wiscons}
\address{Department of Mathematics and Statistics\\
California State University, Sacramento\\
Sacramento, CA 95819, USA}
\email{joshua.wiscons@csus.edu}
%\thanks{The work of the second author was partially supported by the National Science Foundation under grant No.~DMS-1954127.}
\date{\today}

\begin{abstract}
We show that a generically sharply $t$-transitive permutation group of finite Morley rank on a set of rank $r$ satisfies $t\le r+2$ provided the pointwise stabilizer of a generic $(t-1)$-tuple is an $L$-group, which holds, for example, when this stabilizer is solvable or when $r\le 5$. This makes progress towards establishing the natural bound on $t$ implied by the Borovik-Cherlin conjecture that every generically $(r+2)$-transitive permutation group of finite Morley rank on a set of rank $r$ is of the form $\PGL_{r+1}(F)$ acting naturally on $\mathbb{P}^r(F)$. 

Our proof is assembled from three key ingredients that are independent of the main theorem---these address actions of $\Alt(n)$  on  $L$-groups of finite Morley rank, generically $2$-transitive actions with abelian point stabilizers, and simple groups of rank $6$.
\end{abstract}
\maketitle

% % % % % % % % % % % % % % % % % % % %
% % % % % % % % % % % % % % % % % % % %
% SECTION 
% % % % % % % % % % % % % % % % % % % %
% % % % % % % % % % % % % % % % % % % %
\section{Introduction}

The study of generically $t$-transitive actions was popularized by Popov in the context of algebraic groups~\cite{PoV07} and then by Borovik and Cherlin in the more general setting of groups of finite Morley rank~\cite{BoCh08}, although the first substantial results appeared nearly two decades earlier with work of Gropp~\cite{GrU92}. While ordinary $t$-transitivity requires that a group $G$ act transitively on $X^t$ off of the diagonal, 
\emph{generic $t$-transitivity} only requires $G$ act transitively on $X^t$ off of some set of smaller dimension. % (or Morley rank). 

This is considerably more natural than ordinary $t$-transitivity.
For example, the action of $\PGL_{n+1}(F)$ on $\mathbb{P}^{n}(F)$ is generically $(n+2)$-transitive (since the set of projective bases form a single orbit), but it is not even $3$-transitive once $n>1$. 
Generic $t$-transitivity also lends itself to a broader range of applications; this is  particularly well demonstrated by its recent use in the theory of algebraic differential equations \cite{FrMo23}.

% % % % % % % % % % % % % % % % % % % %
% SUBSECTION 
% % % % % % % % % % % % % % % % % % % %
\subsection{The problem}
We work in the setting of groups of finite Morley rank, denoting Morley rank  by $\rk$; the general theory can be found in \cite{BoNe94}. As is common, our use of \emph{definable} includes interpretable; it can simply be read as \emph{algebraic} or \emph{constructible} by those familiar with the algebraic setting but not with Morley rank. For a group acting on a set $X$, we say that $(G,X)$ is a \emph{permutation group of finite Morley rank} if the action is faithful and each of $G$, $X$, and the action of $G$ on $X$ are  definable in some ambient structure of finite Morley rank.

\begin{definition*}
A permutation group $(G, X)$ of finite Morley rank is \emph{generically $t$-transitive} if $G$ has an orbit $\mathcal{O}$ on $X^t$ such that $\rk(X^t \setminus \mathcal{O}) < \rk (X^t)$. The maximum such $t$ for which the action is generically $t$-transitive is called the \emph{degree of generic transitivity} (or \emph{generic transitivity degree} as in~\cite{PoV07}), denoted $\gtd(G, X)$. 
\end{definition*}

%It is a classical result that ordinary $t$-transitivity does not exist for actions on an infinite set once $t\ge 4$~\cite{TiJ52,HaM54}, but 
As we have already observed, generic $t$-transitivity exists very naturally for all $t$; however, one may still look to bound $t$ relative to a fixed $G$ or relative to actions on sets of a fixed rank. 

In the former case, for fixed $G$, it is straightforward to see that $\gtd(G, X)$ is bounded above by $\frac{\rk G}{\rk X}$, but determining the precise value can be nontrivial. One outcome of Popov's work was that, for a fixed connected nonabelian reductive group $G$, the maximum value of $\gtd(G, X)$, taken over all algebraic actions, is in fact associated to a natural action: the action of $G$ on the cosets of some maximal parabolic subgroup. Popov also computes $\gtd(G, X)$ for all connected simple $G$ with $X$ the coset spaces of a maximal parabolic.  

In the case of actions on sets of a fixed rank $r$, it is not immediately clear that there is any bound at all on the degree of generic transitivity, but, as shown by Borovik and Cherlin, there is.  This is essentially the main result of \cite{BoCh08} where they then raise the obvious question of what, precisely, the least upper bound is. In fact, they ask for  more: to identify those actions that achieve the extreme.

\begin{mainproblem}[{\cite[Problem~9]{BoCh08}}]
Assume $(G,X)$ is a transitive and generically $t$-transitive permutation group of finite Morley rank with $X$ of rank $r\ge 1$. Show that if $t\ge r +2$, then $(G,X)$ is isomorphic to $\PGL_{r+1}(F)$ acting naturally on $\mathbb{P}^r(F)$ (and hence $t= r+2$).
\end{mainproblem}

This problem has only been solved in a handful of settings: when $r\le 2$ (by Hrushovski in rank $1$ and the present authors in rank $2$ \cite{AlWi18}), in the theory of algebraically closed fields of characteristic $0$ (by Freitag and Moosa \cite[Theorem~6.3]{FrMo23} leveraging Popov's work), and  in the theory of differentiably closed fields of characteristic $0$ (by Freitag, Jimenez, and Moosa \cite[Theorem~4.3]{FJM23}). It should be mentioned that the \hyperlink{h:problem}{Main Problem} has a natural analog for modules, and this has recently been solved in full generality by Berkman and Borovik \cite{BeBo22}. 

The approach of the authors in rank $2$ (and also of Berkman and Borovik for modules) indicate a canonical case division when attacking the \hyperlink{h:problem}{Main Problem}: is the action generically \emph{sharply} $t$-transitive or not.

\begin{definition*}
A permutation group $(G, X)$ of finite Morley rank is \emph{generically sharply $t$-transitive} if it is generically  $t$-transitive with large orbit $\mathcal{O}$, and $G$ acts without fixing any tuple in  $\mathcal{O}$.
\end{definition*}

The  point is that the natural action of $\PGL_{r+1}(F)$ on $\mathbb{P}^{r}(F)$ is in fact generically \emph{sharply} $(r+2)$-transitive. As such, the ``sharp version'' of the \hyperlink{h:problem}{Main Problem} aims at identification; whereas, the ``non-sharp version'' seeks a contradiction.  The sharp version itself can be broken down further: establish the desired bound on $t$ and then proceed with identification. This  naturally breaks the \hyperlink{h:problem}{Main Problem}  into three sub-problems.

\begin{problem}[Sharpness]\label{prob.Sharpness}
Assume $(G,X)$ is a generically $t$-transitive permutation group of finite Morley rank with $\rk X = r$. Show that if $t\ge r +2$, then the action must be generically \emph{sharply} $t$-transitive.
\end{problem}

\begin{problem}[Bound]\label{prob.BoundSharp}
Assume $(G,X)$ is a generically sharply $t$-transitive permutation group of finite Morley rank with $\rk X = r$. Show that if  $t\ge r +2$, then $t= r +2$.
\end{problem}

\begin{problem}[Identification]\label{prob.Identification}
Assume $(G,X)$ is a transitive and generically sharply $t$-transitive permutation group of finite Morley rank with $\rk X = r$.  Show that if $t = r+2$, then $(G,X)\cong (\PGL_{r+1}(F),\mathbb{P}^r(F))$.
\end{problem}

% % % % % % % % % % % % % % % % % % % %
% SUBSECTION 
% % % % % % % % % % % % % % % % % % % %
\subsection{The theorem}
At present, little has been done around Problem~\ref{prob.Sharpness}, but its solution for modules provides a sliver of hope. Earlier work of the authors addressed aspects of Problem~\ref{prob.Identification} in arbitrary rank \cite{AlWi19}; here, we study Problem~\ref{prob.BoundSharp}, which previously was only handled when $r\le 2$ as consequence of \cite{AlWi18} (see also the  early work of Gropp~\cite{GrU92}). 

Our main result is rank free but comes with a condition on the stabilizer of a generic $(t-1)$-tuple from $X$. In the target case when $(G,X)\cong (\PGL_{r+1}(F),\mathbb{P}^r(F))$, this  stabilizer is precisely a maximal torus of $G$, and for our result, we ``simply'' ask that the stabilizer does not involve a certain (conjecturally nonexistent) class of simple groups. In particular, the condition will automatically be met if this stabilizer is known to be solvable. The precise framework we work in is that of $L$-groups.

\begin{definition*}
Let $H$ be a group of finite Morley rank.
\begin{itemize}
    \item $H$ has \emph{degenerate type} if it has no involutions;
    \item $H$ has \emph{even type} if it contains an infinite elementary abelian $2$-group but no Pr\"ufer $2$-group;
    \item $H$ has \emph{odd type} if it contains a Pr\"ufer $2$-group but no infinite elementary abelian $2$-group;
    \item $H$ has \emph{mixed type} if it contains an infinite elementary abelian $2$-group and a Pr\"ufer $2$-group.
\end{itemize}
We say $H$ is an \emph{$L$-group} if every  definable simple section of $H$ that is of odd type is isomorphic to an algebraic group over an algebraically closed field.
\end{definition*}

It is known---and nontrivial---that every connected group of finite Morley rank has one of the four types  listed above \cite{BoPo90,BBC07} and that the infinite simple ones cannot have mixed type \cite{ABC08}.

\begin{remark}
The \emph{Algebraicity Conjecture} of Cherlin and Zilber posits that every infinite simple group of finite Morley rank is isomorphic to an algebraic group over an algebraically closed field. This analogue of the classification of the finite simple groups remains far from being established, but remarkably it has been confirmed in even and mixed type (with no simple groups of the latter type)~\cite{ABC08}. Attacking degenerate-type simple groups amounts to proving a Feit-Thompson result in this context, which has not seen much progress. The $L$-hypothesis aims to push the theory forward without assuming knowledge of degenerate type groups. %As such, results in this context require that degenerate-type groups be specifically addressed, which, 
In the setting of Problem~\ref{prob.BoundSharp}, this is made possible by \cite{AlWi24}. 
\end{remark}

Here is our main result. 
\begin{maintheorem}
Assume $(G,X)$ is a generically sharply $t$-transitive permutation group of finite Morley rank with $\rk X = r$. 
Let $H$ be the pointwise stabilizer of a generic $(t-1)$-tuple. If $H$ is an $L$-group, then $t\le r +2$.
%If $t\ge r +2$, then $t= r +2$, and $H$ is solvable.
\end{maintheorem}

We stress that the $L$-hypothesis is \emph{not} placed on $G$ but rather only on $H$, which conjecturally should be a torus when $t\ge r +2$. 
As mentioned above, it is  \cite{AlWi24} that provides traction for the $L$-hypothesis; in particular, \cite[Corollary]{AlWi24} tells us at the outset that  either $H$ is solvable or contains involutions. 

\begin{remark}
Our work also yields information when $t= r +2$. In this extremal case,  Propositions~\ref{prop.SymAltActionOnLGroup} and \ref{prop.Rank3HMustBeSolvable} can be used to show that a generic $(t-1)$-tuple must necessarily be solvable except possibly when $r=4$. Analysis similar to what we undertake for $r=3$ in Proposition~\ref{prop.Rank3HMustBeSolvable} (leveraging the significant Proposition~\ref{Prop.SimpleRank6}) may establish solvability when $r=4$, but we do not take it up here. 
\end{remark}

Using the classification of groups of rank at most $5$ \cite{ReJ75,ChG79,WiJ16,DeWi16,FrO18} and the classification of the simple groups of finite Morley rank of even type \cite{ABC08}, the \hyperlink{h:theorem}{Theorem} immediately yields the following.

\begin{maincorollary}
Assume $(G,X)$ is a generically sharply $t$-transitive permutation group of finite Morley rank with $\rk X = r$. 

Then $t\le r +2$ provided any one of the following hold:
\begin{itemize}
    \item $r\le 5$;
    \item $G$ has even type;
    \item the stabilizer of a generic $(t-1)$-tuple is solvable.
\end{itemize}
\end{maincorollary}

% % % % % % % % % % % % % % % % % % % %
% SUBSECTION 
% % % % % % % % % % % % % % % % % % % %
\subsection{The strategy}\label{s.Strategy}
Let us outline our approach to the \hyperlink{h:theorem}{Theorem}. Assume we are in the setting of the \hyperlink{h:theorem}{Theorem} with  $(1,\ldots,t)$ a generic $t$-tuple in $X^t$  and $G_{[t-1]}$ the pointwise stabilizer of $1,\ldots,t-1$. Let $\Sigma$ be the setwise stabilizer of $1,\ldots,t$ and $\Sigma_t \le \Sigma$ the stabilizer of $t$ in $\Sigma$. Generic sharp $t$-transitivity readily implies  $\Sigma\cong\Sym(t)$ and $G_{[t-1]}$ is connected of rank $r$. By  \cite[Lemma~4.27]{AlWi18}, $\Sigma_t$ acts \emph{faithfully} on $G_{[t-1]}$.
Towards a contradiction, we assume that $t\ge r+3$. The main point then is that  $G_{[t-1]}$ carries a faithful action of the relatively large symmetric group $\Sym(r+2)$, and this is  restrictive. 

In Section~\ref{s.Alt(n)OnNonsolvable}, we study actions of $\Alt(n)$ and  $\Sym(n)$ on connected \emph{nonsolvable} $L$-groups of rank $r$ and show in Proposition~\ref{prop.SymAltActionOnLGroup} that typically $n\le r$. When $t\ge r+3$, this yields that either: (1) $G_{[t-1]}$ is  solvable, or (2) we are in the exceptional case of $\Sigma_t \cong \Sym(5)$ and $G_{[t-1]} \cong \PGL_2(F)$ with $\charac F=5$  (which is an expected complication as $\Sym(5) \cong \PGL_2(\mathbb{F}_5)$). 
To push further, both cases require us to work with a larger chunk of $G$ than just $G_{[t-1]}$, but it turns out that $G_{[t-2]}$ (the pointwise stabilizer of $1,\ldots,t-2$) suffices. 

When $G_{[t-1]}$ is solvable, the condition $t\ge r+3$ forces $G_{[t-1]}$ to be an elementary abelian $p$-group. This configuration is studied at the level of $G_{[t-2]}$ and eliminated by Proposition~\ref{prop.Gen2TransAbelianPS}, which looks at generically $2$-transitive actions with abelian point stabilizers. The case of $G_{[t-1]} \cong \PGL_2(F)$ corresponds to $\rk G_{[t-2]} = 6$, and to analyze this, we undertake a rather general study of simple groups of rank $6$ in Section~\ref{s.SimpleRank6}. This culminates with Proposition~\ref{Prop.SimpleRank6}, showing that such groups are so-called $N_\circ^\circ$-groups in the sense of \cite{DeJa16}, which we leverage to obtain a contradiction.

We emphasize that the key ingredients of our proof---Propositions~\ref{prop.SymAltActionOnLGroup}, \ref{prop.Gen2TransAbelianPS}, and~\ref{Prop.SimpleRank6}---live  beyond the specific context of the \hyperlink{h:theorem}{Theorem} and should be interesting in their own right. 

\subsection{Regarding the \texorpdfstring{$L$}{L}-hypothesis}
We would like to highlight that the $L$-hypothesis in the \hyperlink{h:theorem}{Theorem} is only used when we are applying Proposition~\ref{prop.SymAltActionOnLGroup}; or more to the point, removing the hypothesis from Proposition~\ref{prop.SymAltActionOnLGroup} would remove it from the \hyperlink{h:theorem}{Theorem}.

Despite some effort trying to free Proposition~\ref{prop.SymAltActionOnLGroup} from the $L$-hypothesis, we were not able to realize it. What is needed is an extension of Corollary~\ref{SymAltActionOnSimple} from actions of $\Alt(n)$ on simple algebraic groups to actions on simple odd-type groups of finite Morley rank. 
One angle we entertained was to adapt the analysis of actions on degenerate type groups from \cite{AlWi24}; there, the approach is inductive, driven by the fact that a degenerate-type group  acted upon by a Klein $4$-group $K$ is generated by the centralizers of the involutions in $K$~\cite[Theorem~5]{BBC07}. Our attempt to generalize this to actions on odd-type groups failed to find  traction, and this may well be a difficult problem. 

Although the value is unclear, one could try replacing the $L$-hypothesis of Proposition~\ref{prop.SymAltActionOnLGroup} with the \emph{assumption} that the ``four-group generation'' result holds for the simple definable sections of $H$, leading to an analysis of $V$-groups (as $K$ is already taken). We thought very briefly on this. A $V$-assumption moves the analysis of \cite{AlWi24} along quite smoothly until it requires a lower bound on the corank of a proper definable subgroup (see \cite[Corollary~2.5]{AlWi24}), which uses that the degree of generic transitivity of any action of a connected degenerate type group is at most $1$. For odd-type groups there is no such universal bound, but something can still be said. For example, since the degree of generic transitivity of a connected solvable group is at most $2$, one obtains that the degree of generic transitivity of a connected minimal simple group is at most $3$. This may be enough to prove a statement along the following lines: if $\Alt(n)$, with $n$ sufficiently large, acts on a connected minimal simple $V$-group $H$ of finite Morley rank, then $n\le \rk H + 2$. 

There is also the question of trying to extend  Corollary~\ref{SymAltActionOnSimple} to simple $L^*$-groups $G$ of odd type. The $L^*$-hypotheses assumes all \emph{proper} definable sections of $G$ are $L$-groups but makes no claim on $G$ itself. The  case of simple $L^*$-groups with sufficiently high Pr\"ufer rank has recently been reviewed and clarified by Cherlin: assuming such a $G$ has Pr\"ufer $2$-rank  at least $3$ and that there are certain internal restrictions on how $2$-tori act (the \textit{NTA}$_2$ hypothesis), then either $G$ is algebraic or contains a strongly embedded subgroup~\cite{ChG24}.  As such, a test case for extending Corollary~\ref{SymAltActionOnSimple} to the $L^*$-setting might be when $G$ has high Pr\"ufer $2$-rank and a strongly embedded subgroup. Or, instead focusing on low Pr\"ufer $2$-rank, one could try to analyze actions on 
groups of type $\cibo_1$, $\cibo_2$, or $\cibo_3$ (as defined in Section~\ref{s.SimpleRank6}) with an initial focus on the minimal simple case.

% % % % % % % % % % % % % % % % % % % %
% % % % % % % % % % % % % % % % % % % %
% SECTION 
% % % % % % % % % % % % % % % % % % % %
% % % % % % % % % % % % % % % % % % % %
\section{Three key ingredients}

Here we develop the main ingredients of our proof of the \hyperlink{h:theorem}{Theorem}, as laid out in Section~\ref{s.Strategy}. These results are independent of the \hyperlink{h:theorem}{Theorem} and address: actions of $\Alt(n)$ and  $\Sym(n)$ on  $L$-groups of finite Morley rank (Section~\ref{s.Alt(n)OnNonsolvable}), generically $2$-transitive actions with abelian point stabilizers (Section~\ref{s.PermGroups}), and simple groups of rank $6$ (Section~\ref{s.SimpleRank6}).

% % % % % % % % % % % % % % % % % % % %
% SUBSECTION 
% % % % % % % % % % % % % % % % % % % %
\subsection{Actions of \texorpdfstring{$\Alt(n)$}{Alt(n)} on nonsolvable groups}\label{s.Alt(n)OnNonsolvable}

Our work here focuses the following conjecture, which previously had only been addressed for actions on groups of degenerate type~\cite{AlWi24}. 

\begin{conjecture*}[see {\cite[Section~1.2]{CDW23}}]
Suppose that $\Alt(n)$ acts definably and faithfully by automorphisms on a nonsolvable connected group  $H$ of finite Morley rank. Then $n \le \rk H$ for sufficiently large  $n$.
\end{conjecture*}

The solvable version of the conjecture allows for $n \le \rk H + 2$ and was addressed in \cite[Lemma~2.7]{AlWi24} with the critical case of abelian $H$ taken care of by the analysis of \cite{CDW23}.

Leveraging the algebraic theory, we solve the conjecture for actions on $L$-groups; this is Proposition~\ref{prop.SymAltActionOnLGroup}. Our proof consists to a large degree of simply collecting and organizing known results. As indicated in Corollary~\ref{cor.ngeninefamily}, the bound of $\rk H$ on $n$ is typically rather weak and should be more on the order of $\sqrt{\rk H}$ for generic $n$.

\begin{proposition}\label{prop.SymAltActionOnLGroup}
Assume $\Alt(n)$ or $\Sym(n)$ acts definably and faithfully by automorphisms on a nonsolvable connected $L$-group  $H$ of finite Morley rank.

 If $n\ge 9$, then $n\le \rk H$. If the acting group is $\Sym(n)$, then  $n\ge 5$ also implies $n\le \rk H$ unless $n=5$ and $H\cong Z(H) \times \PGL_2(F)$ with $\charac F = 5$. 
\end{proposition}

Our approach will consider projective representations of $\Alt(n)$ of minimal dimension for which we first collect results on the linear representations of $\Alt(n)$ for $n\ge 8$ and its unique Schur cover, which we denote by $\widehat\Alt(n)$.  

\begin{fact}[{\cite{WaA76,WaA77}}]
Let $n\ge 9$. Suppose $\Alt(n) \le \GL_d(F)$ for $F$ an algebraically closed field. Then $d\ge n-1-\kappa_n$ where $\kappa_n = 1$ if $\charac F$ divides $n$ and $\kappa_n = 0$ otherwise. Moreover, the same bound holds for $n=7,8$ provided $\charac F \neq 2$.
\end{fact}

\begin{fact}[{\cite[Theorem A]{KlTi04}}]
Let $n\ge 8$. Suppose $\widehat\Alt(n) \le \GL_d(F)$ for $F$ an algebraically closed field of characteristic not $2$. Then $d \ge 2^{\lfloor(n-2-\kappa_n)/2\rfloor}$ where $\kappa_n = 1$ if $\charac F$ divides $n$ and $\kappa_n = 0$ otherwise.
\end{fact}

The previous facts combine to yield the following. 

\begin{corollary}\label{cor.EmbedAlt(n)InPGL}
Let $n\ge 9$. Suppose $\Alt(n) \le \PGL_d(F)$ for $F$ an algebraically closed field. Then $n\le d+2$.
\end{corollary}

We now consider embeddings of $\Alt(n)$ into other simple algebraic groups. Our approach is cheap: given a simple algebraic group $G$, we embed $G$ into $\PGL_d(F)$ and then apply Corollary~\ref{cor.EmbedAlt(n)InPGL}.

\begin{fact}\label{fact.EmbedSimpleAlgebraicInPGL}
Suppose $G$ is a (group-theoretically) simple algebraic group over an algebraically closed field $F$. Then Table~\ref{Table:MinimalRepresentations} gives the algebraic dimension of $G$ and minimal $d$ such that  $G$ embeds into $\PGL_d(F)$.
\end{fact}
\begin{table}[ht] 
\small
\centering
\renewcommand{\arraystretch}{1.25}
    \begin{tabular}{c|c|c}
       Type of $G$  & $\dim(G)$ & $d$ such that $G$ embeds in $\PGL_d(F)$ \\ \hline
       $A(\ell)$  & $(\ell+1)^2 - 1$ & $\ell+1$ \\ \hline
       $B(\ell)$  & $2\ell^2 + \ell$ & $2\ell+1$ \\ \hline
       $C(\ell)$  & $2\ell^2 + \ell$ & $2\ell$ \\ \hline
       $D(\ell)$  & $2\ell^2 - \ell$ & $2\ell$ \\ \hline
       $E(6)$  & 78 & $27$ \\ \hline
       $E(7)$  & 133 & $56$\\ \hline
       $E(8)$  & 248 & 248 \\ \hline
       $F(4)$  & 52 & $26$\\ \hline
       $G(2)$  & 14 & $7$ \\ \hline
    \end{tabular}
\caption{Dimensions of the simple algebraic groups and minimal $d$ such that they embed into $\PGL_d(F)$}
\label{Table:MinimalRepresentations}
\end{table}

A few remarks are in order to justify the data contained in Table~\ref{Table:MinimalRepresentations}. To obtain the second column, it suffices to compute the cardinality of the set of positive roots in the corresponding linear algebraic group. The dimension of the group is then twice this cardinality plus the dimension of a maximal torus. The numerical data for this can be found in Section 12.2 of \cite{HuJ72} (see also Section 28 of \cite{HuJ75}). A compact exposition of the dimensions in the third column, in a more general setting, can be found in \cite{Che98} with the justification given in Theorem 1.2 of that paper. (It seems that a slight typo has found its way into the formula for $F_4$ on page 672 of that paper, but the argument on page 676 rectifies it.)

\begin{corollary}\label{cor.ngeninefamily}
Let $n\ge 9$. Suppose $\Alt(n) \le G$ for $G$ an $r$-dimensional, simple algebraic group over an algebraically closed field. 
Then $n\leq r$, and if $G$ is a classical group, then $n\le f(r)$ where $f(r)$ is given in Table~\ref{Table:BoundingN}. 
\end{corollary}
\begin{table}[ht]
\small
\centering
\renewcommand{\arraystretch}{1.5}
    \begin{tabular}{c|c}
       Type of $G$  & $f(r)$ \\ \hline
       $A(\ell)$  &  $2+ \sqrt{r+1}$\\ \hline
       $B(\ell)$  &  $\frac{5+ \sqrt{8r+1}}{2}$\\ \hline
       $C(\ell)$  &  $\frac{3+ \sqrt{8r+1}}{2}$\\ \hline
       $D(\ell)$  &  $\frac{5+ \sqrt{8r+1}}{2}$\\ \hline
    \end{tabular}
\caption{Upper bounds on $n$ when $\Alt(n)$ embeds into $r$-dimensional classical groups}
\label{Table:BoundingN}
\end{table}

% \begin{corollary}\label{cor.ngenineexceptional}
% Let $n\ge 9$. Suppose $\Alt(n) \le G$ for $G$ a simple exceptional algebraic group over an algebraically closed field. Then $n\leq \dim G$.
% \end{corollary}
\begin{proof}
Corollary~\ref{cor.EmbedAlt(n)InPGL} and Fact~\ref{fact.EmbedSimpleAlgebraicInPGL}  resolve this for all cases except $E(8)$, so suppose $G$ is of type $E(8)$. The Weyl group involves only the primes up to $7$, so if $E<G$ is an elementary abelian $p$-group with $p>7$ and $p\neq \charac F$, then $E$ is toral, implying that the $p$-rank of such an $E$ is at most $8$. 

If $n\ge 99$, then $\Alt(n)$ contains an elementary abelian $11$-group of $11$-rank $9$, showing that  $n< 99$ when $\charac F \neq 11$. And if $\charac F = 11$, we may instead consider elementary abelian $13$-groups to see that $n< 117$. 

The proof is complete, but we note that this method can be applied to the primes $2$ and $3$ (to produce much smaller bounds) provided one also accounts for nontoral such $E$, for which \cite{GrR91} gives the necessary information.
\end{proof}

\begin{corollary}\label{SymAltActionOnSimple}
Assume $\Alt(n)$ or $\Sym(n)$ acts definably and faithfully on  a simple algebraic group $G$ over an algebraically closed field.

If $n\ge 9$, then $n\le \dim G$. If the acting group is $\Sym(n)$, then  $n\ge 5$ also implies $n\le \dim G$ unless $n=5$ and $G=\PGL_2(F)$ with $\charac F = 5$. 
\end{corollary}
\begin{proof}
By \cite[Ch.~II,~Fact~2.25]{ABC08} (generalizing \cite[Theorem~27.4]{HuJ75}), the acting group is contained in $\Inn(G) \cdot \Gamma$ for $\Gamma$ the group of graph automorphisms relative to a fixed maximal torus and Borel subgroup, so  $\Alt(n)$ embeds into $G$ in the cases under consideration. By Corollary~\ref{cor.ngeninefamily}, it only remains to address when $5\le n \le 8$ and the acting group is $\Sym(n)$. Let $r = \dim G$. If $r\ge 8$, there is nothing to show, so we may assume $G$ is of type $A(1)$. In particular, $r= 3$, and $G=\PGL_2(F)$ has no graph automorphisms. Thus, $\Sym(n)$ embeds into $\PGL_2(F)$. This is now quite classical. Via the adjoint representation of $\PGL_2(F)$, we have a faithful $3$-dimensional representation of $\Sym(n)$, so the only possibility is $n=5 = \charac F$  (see for example~\cite{DiL08}).
\end{proof}

\begin{remark}
Since the Morley rank of an algebraic group $G$ is always greater than or equal to its algebraic dimension, the bound $n\le \dim G$ appearing in  Corollary~\ref{SymAltActionOnSimple} can be replaced by $n\le \rk G$.
\end{remark}

\begin{proof}[Proof of Proposition~\ref{prop.SymAltActionOnLGroup}]
We will denote by $X$ the acting group, $n$ the degree of the permutation group $X$, and $r$ the rank of $H$. Also, set $X' = [X,X] = \Alt(n)$, and let $\sigma(H)$ denote the solvable radical of $H$.

We begin by addressing the exceptional case of $n=5$.

\begin{claim}\label{claim.prop.SymAltActionOnLGroup.n=5}
If $X= \Sym(5)$ and $r<5$, then $H\cong Z(H) \times \PGL_2(F)$ with $F$ of characteristic $5$. 
\end{claim}
\begin{proofclaim}
Since $H$ is nonsolvable and connected groups of Morley rank at most $2$ are
solvable~\cite{ChG79}, $r$ is either $3$ or $4$.

First assume $r=3$. Then $H$ is quasisimple with a finite center. Notice that $X$ acts faithfully on $H/Z(H)$;  otherwise,  the connected group $[X',H]$ is contained in $Z(H)$ and hence is trivial (since $Z(H)$ is finite), contradicting faithfulness of the action on $H$.
Now, simple groups of rank $3$ are known to be algebraic (\cite{ChG79,FrO18}), so Corollary~\ref{SymAltActionOnSimple} applies to $H/Z(H)$. Thus, $H/Z(H)\cong\PGL_2(F)$ with $\charac F = 5$. Moreover, by the theory of central extensions (see \cite[Ch.~II,~Lemma~2.21]{ABC08}), $H$ is algebraic, hence of the form $\PGL_2(F)$ or $\SL_2(F)$. In either case, $X$ must act by inner automorphisms (see \cite[Ch.~II,~Lemma~2.27]{ABC08}), but $\SL_2(F)$ does not embed $\Sym(5)$. Thus, $Z(H) = 1$, and $H$ has the desired structure.

Next consider when $r=4$. By \cite{AlWi24}, $H$ contains involutions, and its  structure is given by  \cite{WiJ16}: $H = Z^\circ(H)*Q$ with $Q$ quasisimple of rank $3$. Since $Z^\circ(H)$ has rank $1$, it is centralized by $X'$, so $X$ must act faithfully on $Q$. Our analysis above then shows $Q \cong\PGL_2(F)$ in characteristic $5$, which in turn forces  $H = Z(H)\times Q$.
\end{proofclaim}

For the remainder of the proof, we assume $n\ge 6$ and proceed by contradiction. Thus we assume $r<n$, and we choose $H$ to be a counterexample to Proposition~\ref{prop.SymAltActionOnLGroup} of minimal rank $r$. 

\begin{claim}
$\sigma^\circ(H)=1$.
\end{claim}
\begin{proofclaim}
Assume $\sigma^\circ(H)\neq1$. Then $\rk(H/\sigma^\circ(H))<r$, and since $X$ acts  on $H/\sigma^\circ(H)$,  minimality of $r$ implies that this action is not faithful. Thus $[X',H]\leq\sigma^\circ(H)$.
Now, if we also have $[X',\sigma^\circ(H)]=1$, then $[X',[X',H]] = 1$, and the three subgroups lemma implies that the group $[X',X'] = X'$ centralizes $H$, against faithfulness of the action on $H$. Thus, $[X',\sigma^\circ(H)]\neq 1$, so $X$ acts faithfully on $\sigma^\circ(H)$. 

Then \cite[Lemma~2.7]{AlWi24} implies $\rk(\sigma^\circ(H))\geq n-2$ in all cases except possibly when $X=\Alt(9)$. However, if $X=\Alt(9)$ and $\rk(\sigma^\circ(H)) < 7$, we can use that the minimal rank of a faithful connected $\Alt(9)$-module is $7$ (see the main result of \cite{CDW23} and the final remark of Section~3.3, also in \cite{CDW23}) to eliminate this case. Thus $\rk(\sigma^\circ(H))\geq n-2$, meaning that $\rk(H/\sigma^\circ(H))\le 2$. Since all connected groups of Morley rank at most $2$ are solvable, this is a contradiction, so $\sigma^\circ(H)=1$.
\end{proofclaim}

\begin{claim}
We may assume $H$ is semisimple: there exist $k\geq1$ and connected simple groups of finite Morley rank
$H_1,\ldots,H_k$ such that $H=H_1\times\cdots\times H_k$.
\end{claim}
\begin{proofclaim}
Since $\sigma^\circ(H)=1$, we find (as in the proof of Claim~\ref{claim.prop.SymAltActionOnLGroup.n=5}) that $X$ acts faithfully on $H/\sigma(H)$. Thus $H/\sigma(H)$ is also a counterexample of rank $r$, and the socle of $H$ is as well. Thus we may assume $H$ equals its socle, which is semisimple by \cite[Theorem~7.8]{BoNe94}.
\end{proofclaim}

\begin{claim}\label{claim.prop.SymAltActionOnLGroup.n>=9}
If $n\ge 9$, we have a contradiction.
\end{claim}
\begin{proofclaim}
Since $H$ is an $L$-group, each $H_i$ is either a simple linear algebraic group over an algebraically closed field or of degenerate type. By the inductive assumption, $n>r\geq k$ (in fact, $n>r\geq 3k$ since 
the $H_i$ are not solvable). As $X'$ is simple and permutes the set $\{1,\ldots,k\}$, the only possible action is the trivial one. 
Thus $X'$ normalizes each $H_i$. However, $\rk H_i \le r <n$, so 
$X'$ centralizes each $H_i$ using Corollary~\ref{SymAltActionOnSimple} in the algebraic case and \cite{AlWi24} in the degenerate case. This contradicts faithfulness. 
\end{proofclaim}

\begin{claim}
If $6\le n\le 8$ and $X= \Sym(n)$, we have a contradiction.
\end{claim}
\begin{proofclaim}
By the contradictory assumption, $r\leq 7$. Since connected groups of Morley rank at most $2$ are
solvable, $k$ is at most $2$. 

We first eliminate the case when $k=2$. If $k=2$ then one of the groups, say $H_1$, is of rank $3$
and the other is of rank $4$. As a result $X$ cannot permute these two groups so normalizes each one. Minimality of
$r$ implies that the action of $X$ on both components has a nontrivial kernel, which must be $X'$. It follows that the action of $X$ on $H$ is not faithful, a contradiction. 

Thus, $k=1$ and $H$ is simple.
As in Claim~\ref{claim.prop.SymAltActionOnLGroup.n>=9}, Corollary~\ref{SymAltActionOnSimple} and \cite{AlWi24} provide a contradiction.
\end{proofclaim}

The proof of Proposition~\ref{prop.SymAltActionOnLGroup} is now complete.
\end{proof}
\setcounter{claim}{0} 

% % % % % % % % % % % % % % % % % % % %
% SUBSECTION 
% % % % % % % % % % % % % % % % % % % %
\subsection{Actions with virtually-abelian point-stabilizers}\label{s.PermGroups}

\begin{proposition}\label{prop.Gen2TransAbelianPS}
Let $(G,X)$ be a transitive and generically $2$-transitive permutation group of finite Morley rank with $G$ connected. 

If $G_x^\circ$ is abelian, then every virtually definably primitive quotient $\overline{X}$, with corresponding kernel $N$, satisfies $(G/N,\overline{X})\cong (\AGL_1(L),L_+)$ and $G_xN/N \cong L^\times$ for $L$ an algebraically closed field.
\end{proposition}

\begin{setup}
Throughout this section we assume $(G,X)$ is a transitive permutation group of finite Morley rank.
\end{setup}

\begin{lemma}\label{lem.GenNTransAbelianPS}
Let $(G,X)$ be generically $t$-transitive with $t\ge 2$. Assume $(G_{1,\ldots,t-1})^\circ$ is abelian for $(1,\ldots,t)$ in the generic orbit of $G$ on $X^t$. 

Then every definable quotient $\overline{X}$, with corresponding kernel $N$, satisfies  $\rk (G_{1,\ldots,{t-1}}\cap N) = \rk X - \rk \overline{X}$. In particular, $\rk X > \rk \overline{X}$ implies $\rk N > \rk X - \rk \overline{X}$.
\end{lemma}
\begin{proof}
Let $\overline{X}$ and $N$ be as given; then $X$ and  $\overline{X}$ are connected by \cite[Lemma~1.8]{BoCh08}. 
Set $H:= (G_{1,\ldots,{t-1}})^\circ$. Since $H$ is abelian and  generically transitive on $X$, $H$ is generically regular (using \cite[Lemma~1.6]{BoCh08}), so $\rk H = \rk X$. By \cite[Lemma~6.1]{BoCh08}, the same is true on $\overline{X}$ modulo the kernel, so $\rk HN/N =  \rk \overline{X}$. Thus, $\rk H - \rk H\cap N = \rk \overline{X}$, and all together we have $\rk H\cap N = \rk H - \rk \overline{X} = \rk X - \rk \overline{X}$.

For the final point, if $\rk X > \rk \overline{X}$, then our work shows $N^\circ \neq 1$, so faithfulness implies that $N^\circ\nsubseteq H$. Thus,  $\rk N > \rk H\cap N$.
\end{proof}

\begin{lemma}\label{lem.Gen2TransAbelianPS}
Let $(G,X)$ be generically $2$-transitive. If $G_x$ is abelian-by-finite, then every  definable quotient $\overline{X}$, with corresponding kernel $N$, satisfies $G^\circ_{\bar{x}} = G^\circ_xN$, so $G_{\bar{x}}/N$ is also abelian-by-finite.
\end{lemma}
\begin{proof}
Let $\overline{X}$ and $N$ be as given. We proceed by induction on the difference $d = \rk X - \rk \overline{X}$.  
If $d=0$, the index of $G_{x}$ in $G_{\bar{x}}$ is finite, and we are done. 

Assume  $d\ge 1$.
Then by Lemma~\ref{lem.GenNTransAbelianPS}, $N^\circ$ is nontrivial, and as $N^\circ$ is not contained in $\rk G_{x}$ (by faithfulness), $\rk G_{x}N > \rk G_{x}$.
Now let $\widetilde{X}$ be the  quotient of $X$ determined by $G_{x}N$ and $K$  the corresponding kernel. We simultaneously study the actions of $G$ on $X$, $\widetilde{X}$, and $\overline{X}$, which correspond to the sequence of stabilizers $G_x <  G_{\tilde{x}} \le G_{\bar{x}}$. We  show that the hypotheses on $X$ pass to $\widetilde{X}$, to which induction applies to obtain the  result for $\overline{X}$.

Since $N \le G_{x}N = G_{\tilde{x}}$,  normality of $N$ forces $N\le K$, and as  $K \le G_{\tilde{x}} \le G_{\bar{x}}$,  we also have $K\le N$, hence equality. Thus $G_{\tilde{x}}/K = G_{\tilde{x}}/N = G_{x}N/N$ is abelian-by-finite, so $(G/K,\widetilde{X})$ satisfies the hypotheses of the lemma.  Induction applies, and as $\overline{X}$ is a quotient of $\widetilde{X}$, $G^\circ_{\bar{x}}/K = G^\circ_{\tilde{x}}N/K$. Using that $K = N$ and $ G_{\tilde{x}} = G_{x}N$, we obtain the desired result.
\end{proof}

\begin{proof}[Proof of Proposition~\ref{prop.Gen2TransAbelianPS}]
Taking $\overline{X}$ and $N$ as stated, Lemma~\ref{lem.Gen2TransAbelianPS} shows that the point stabilizers of $(G/N,\overline{X})$ are again abelian-by-finite, so by \cite[Lemma~3.7]{WiJ16}, they are in fact connected and abelian. Then \cite[Proposition~3.8]{WiJ16} applies to give the desired structure of $(G/N,\overline{X})$ and hence the structure of $G_xN/N$ as well since $G_{\bar{x}} = G^\circ_xN$ (from Lemma~\ref{lem.Gen2TransAbelianPS}).
\end{proof}

Though we do not make use of it here, we record a corollary of Proposition~\ref{prop.Gen2TransAbelianPS} that appears quite relevant to continued work on the \hyperlink{h:problem}{Main Problem}.

\begin{corollary}\label{cor.Gen2TransToralPS}
Let $(G,X)$ be generically $2$-transitive with $G$ connected. If $G_x^\circ$ is a good torus, then $G$ has a definable subnormal series $G = N_0 \trianglerighteq N_1 \trianglerighteq \cdots \trianglerighteq N_m = 1$ such that $N_{i}/N_{i+1} \cong \AGL_1(L_i)$ with each $L_i$ an algebraically closed field.
\end{corollary}
\begin{proof}
Fix $x,y$ in general position.  

We begin with a few observations. First, \cite[Lemma~3.7]{WiJ16} implies $G_x$ is connected. It also implies $G_x$ is a \textit{maximal} good torus, which we now show. Indeed, if not, $G_x$ and $G_y$ would be properly contained in good tori $T_x$ and $T_y$, which must intersect nontrivially by rank considerations. Further, $T_x$ and $T_y$ would generate $G$ since $G_x$ and $G_y$ do (see for example \cite[Lemma~4.12]{AlWi18}), and this would force $G$ to have a nontrivial center, against  \cite[Lemma~3.7]{WiJ16}. For the same reason,  $G_x\cap G_y = 1$, so the action is generically \emph{sharply} $2$-transitive.

We  apply Proposition~\ref{prop.Gen2TransAbelianPS} to the action of $G$ on $X$. Letting $N$ be the kernel of a virtually definably primitive quotient $\overline{X}$, we have $G/N\cong \AGL_1(L)$. Since Proposition~\ref{prop.Gen2TransAbelianPS} applies to \textit{any} virtually definably primitive quotient, we may assume $N$ is connected (by passing from the quotient determined by $G_xN$ to the one determined by $G_xN^\circ$).

Now let $\mathcal{O}$ be the orbit of $N$ on $X$ that contains $x$. We show that Proposition~\ref{prop.Gen2TransAbelianPS} applies to the action of $N$ on $\mathcal{O}$ and then conclude by induction. Only faithfulness and generic $2$-transitivity need to be verified. 

Now, $N_x$ and $N_y$ must be maximal good tori of $N$, so, by conjugacy of maximal tori in $N$, $N_y = N_{y'}$ for some $y' \in \mathcal{O}$. This was the key point. Since $G_x\cap G_y = 1$, we find that $N_x\cap N_{y'} = 1$, implying, in particular, that the action of $N$ on $\mathcal{O}$ is faithful. 
This also can be used  to show that 
the action of $N$ on $\mathcal{O}$ is generically $2$-transitive. Let $r = \rk X$ and $s=\rk \overline{X}$. We have that $r-s = \rk G_{\bar{x}}/G_x = \rk G_xN/G_x = \rk N/N_x = \rk \mathcal{O}$; similarly, we find that $\rk N_x  = r-s = \rk \mathcal{O}$. Thus, since $N_x\cap N_{y'} = 1$, the orbit of $N_x$ containing $y'$ is generic, establishing generic $2$-transitivity of $(N,\mathcal{O})$.
\end{proof}

% % % % % % % % % % % % % % % % % % % %
% SUBSUBSECTION 
% % % % % % % % % % % % % % % % % % % %
\subsection{Groups of rank 6}\label{s.SimpleRank6}

Our analysis in the proof of the \hyperlink{h:theorem}{Theorem} leads to an exceptional case for which we need rather detailed information about quasisimple groups of Morley rank $6$. We prove a reasonably strong result for simple groups of Morley rank $6$ that places the quasisimple groups within the $N_\circ^\circ$-framework of Deloro and Jaligot. 

\begin{definition*}
A group $G$ of finite Morley rank is called an \textbf{$N_\circ^\circ$-group} if $N^\circ_G(A)$ remains \textit{solvable} whenever $A$ is an infinite, definable, connected, and \textit{solvable} subgroup of $G$.
\end{definition*}

\begin{proposition}\label{Prop.SimpleRank6}
    If $G$ is a simple group of Morley rank $6$, then $G$ has no rank $5$ subgroups and every connected rank $4$ subgroup is solvable. Consequently, every quasisimple group of Morley rank $6$ with finite center is an $N_\circ^\circ$-group.
\end{proposition}

In reading Proposition~\ref{Prop.SimpleRank6}, bear in mind that, conjecturally, the only simple groups of Morley rank $6$ are those of the form $\PSL_2(L)$ with $L$ of rank $2$. Also, we note that our use of the proposition occurs in a restricted setting where we will be able to show that connected centralizers of involutions are solvable; this will allow us to bring to bear the main result of \cite{DeJa16}.

\begin{fact}[{\cite{DeJa16}}]\label{fact.DJ16}
Let $G$ be an infinite connected nonsolvable $N_\circ^\circ$-group of finite Morley rank of odd type. Further assume that $C^\circ(i)$ is solvable for all $i\in I(G)$. Then one of the following holds.
\begin{description}
\item[\cibo] $C^\circ(i)$ is a Borel subgroup of $G$, and either
\begin{description}
\item[\cibo$_1$] $\mr_2 (G) = 1$, $C(i)$ is a self-normalizing Borel subgroup of $G$;
\item[\cibo$_2$] $\mr_2 (G)= 2$, $\pr_2 (G) = 1$,  $C^\circ(i)$ is an abelian Borel subgroup of $G$ inverted by any $\omega\in C(i)\setminus\{i\}$, and $\rk G = 3\cdot\rk C^\circ(i)$; or
\item[\cibo$_3$] $\mr_2 (G)= \pr_2 (G) = 2$, $C(i)$ is a self-normalizing Borel subgroup of $G$.
\end{description}
\item[Algebraic] $G \cong \PSL_2(K)$.
\end{description}
\end{fact}

\begin{proof}[Proof of Proposition~\ref{Prop.SimpleRank6}]
Assume $G$ is a simple group of Morley rank $6$. That $G$ has no rank $5$ subgroups is due to Hrushovski, see \cite[Theorem~11.98]{BoNe94}. So we aim to show that connected rank $4$ subgroups are solvable. Once done, the final sentence of Proposition~\ref{Prop.SimpleRank6} then follows from the fact that connected groups of rank at most $2$ are solvable.

Now suppose $M<G$ is nonsolvable and connected of rank $4$. Our general strategy is as follows: (1) use the presence of the nonsolvable $M$ of rank $4$ to produce a well-structured solvable  $D$ of rank $4$, (2) study $(G,G/D)$ to identify $G$ as $\PSL_2(F)$, and (3) note that $\PSL_2(F)$ has no subgroup like $M$.  

\begin{claim}\label{Claim.Prop.SimpleRank6.3Trans}
If $H< G$ is connected (solvable or not) of rank $4$, then $(G,G/H)$ is generically $3$-transitive and $H$ is nonnilpotent.
\end{claim}
\begin{proofclaim}
Suppose $H<G$ with $\rk H = 4$. Then $H$ is a maximal connected subgroup, so \cite[Proposition~2.3]{BoCh08} implies $s\cdot \gtd(G,G/H) \le \rk G \le s\cdot \gtd(G,G/H) + s(s-1)/2$ for $s = \rk G/H$. Thus, $2\cdot \gtd(G,G/H) \le 6 \le 2\cdot \gtd(G,G/H) + 1$, implying that $\gtd(G,G/H) = 3$.
Consequently, $\gtd(H,G/H) = 2$, so \cite[Proposition 4.24]{AlWi18} shows  $H$ is nonnilpotent.
\end{proofclaim}

The claim implies that $(M,G/M)$ is generically $2$-transitive, so our chosen $M$  contains involutions. Its structure is then given by \cite[Corollary~A]{WiJ16} together with \cite{FrO18}: $M = Q*Z$ with $Q\cong\pSL_2(K)$ and $Z = Z^\circ(M)$ (and $\rk Z = 1$).

\begin{claim}\label{Claim.Prop.SimpleRank6.StructureM}
Let $M^g$ be a generic conjugate of $M$. Set $A:=(M\cap M^g)^\circ$ and $C:=C^\circ_G(A)$.
Then
\begin{enumerate}
    \item $\rk A =\rk C = 2$ and $\rk (A\cap C) = 0$;
    \item $C = \langle Z,Z^g \rangle $;
    \item $A/Z(A) \cong \AGL_1(K)$ and $C/Z(C) \cong \AGL_1(L)$ for  fields $K$ and $L$;
    \item $Z \cong L^\times \cong Z^g$.
\end{enumerate}
\end{claim}
\begin{proofclaim}
We first claim that neither $Z$ nor $Z^g$ are contained in $A$. To see this, first note that $A$ is the connected component of a generic $2$-point stabilizer in the action  of $G$ on $G/M$. Thus, if $Z\le A$, then $Z$ fixes a point in the generic orbit of $M$ on $G/M$, so as $Z$ is central in $M$, this implies  $Z$ fixes a generic subset of $G/M$, hence is in the kernel of the action, a contradiction to simplicity of $G$. We conclude that $Z\not\le A$ and similarly that $Z^g\not\le A$.

We next show that $A/Z(A) \cong \AGL_1(K)$. First, observe that $\rk A=2$; this follows from the fact that $A$ is the connected component of a generic $2$-point stabilizer in $(G,G/M)$, which is generically $3$-transitive. Now, since $\rk Z = 1$ and $Z\not\le A$, we have that $\rk (Z\cap A) = 0$, so  
the image of $A$ in $M/Z$ has rank $2$.  By the structure of $M$, $A/Z(A)$ must be isomorphic to a Borel of $\PSL_2(K)$, so $A/Z(A) \cong \AGL_1(K)$, as desired.

We now look at $C$. By the structure of $M$, we certainly have $C\ge \langle Z,Z^g\rangle$, and as we now also know the structure of $A$, we find that  $\rk (C\cap A) = 0$. If $\rk C >2$, then the group $AC$, would have rank at least five. However, $G$ has no rank $5$ subgroups, and $G=AC$ is impossible by simplicity. Thus, $\rk C =2$ and $C= \langle Z,Z^g\rangle$. 

We have observed that $Z\not\le A$, so $Z\not\le M^g$. By the structure of $M$ and rank considerations, $M^g = C^\circ_G(Z^g)$, so $[Z,Z^g] \neq 1$. Since $C$ has rank $2$ and distinct noncommuting rank $1$ subgroups, it must be that $C$ is nonnilpotent, so $C/Z(C) \cong \AGL_1(L)$ for some field $L$. This also shows that $Z,Z^g$ are isomorphic to tori in $\AGL_1(L)$.
\end{proofclaim}

Let $D = A*C$ for $A$ and $C$ as in Claim~\ref{Claim.Prop.SimpleRank6.StructureM}. Set $U_A = F^\circ(A)$, and let $T_A$ be a maximal torus of $A$. Similarly define $U_C$ and $T_C$. Then $A = U_A \rtimes T_A$ with $Z(A) \le T_A$ and similarly for $C$ (see \cite[Theorem~2]{ChG79}).

\begin{claim}\label{Claim.Prop.SimpleRank6.StructureD}
$D = F^\circ(D)\rtimes T$ where $F^\circ(D) = U_AU_C$ and $T = T_AT_C$.
\end{claim}
\begin{proofclaim}
We need only show $F^\circ(D)\cap T = 1$. Let $z\in F^\circ(D)\cap T\le Z(D)$. Writing $z=ac$ with $a\in U_A$ and $c\in U_C$, we find $a\in Z(A)\le T_A$ and $c\in Z(C)\le T_C$. But $U_A \cap T_A = 1 = U_C \cap T_C$, so $z = 1$.
\end{proofclaim}

Set $\hat D = N_G(D)$ and  $Y = G/\hat D$. As $G$ has no rank $5$ subgroups, $\hat D$ is maximal in $G$.
 By Claim~\ref{Claim.Prop.SimpleRank6.3Trans}, $(G,Y)$ is generically $3$-transitive. Let $1,2,3\in Y$ be in general position, with $1$ representing the trivial coset $\hat D$ (so $G_1=\hat D$).

\begin{claim}\label{Claim.Prop.SimpleRank6.2TransOnG/D}
$(G,Y)$ is $2$-transitive. Consequently, $\rk G_{1,2} = 2$.
\end{claim}
\begin{proofclaim}
Since $Y$ is connected of rank $2$, we can show the action is $2$-transitive by showing that $G_1$ has a  unique orbit of rank $0$ (namely $\{1\}$) and no orbits of rank $1$.

As $G_1$ is a maximal subgroup,  the action is primitive, so $\{1\}$ is indeed the unique rank $0$ orbit of $G_1$ (see for example \cite[Lemma~3.4]{WiJ16}).

To show  $G_1$ has no orbits of rank $1$, it suffices to show it for $G^\circ_1$. 
Assume that $G^\circ_1$ has an orbit $\mathcal{O}$ of rank $1$; it also has degree $1$. Choose $y\in \mathcal{O}$, and let $K$ be the kernel of the action of $G^\circ_1$ on $\mathcal{O}$. 
Note that $G^\circ_1$ is generated by its maximal decent tori (using Claim~\ref{Claim.Prop.SimpleRank6.StructureD}), which are conjugate in $G^\circ_1$, so $K$ does not contain a maximal decent torus of $G^\circ_1$. Since $G$ has no rank $5$ subgroups, $\langle G_1^\circ,G_y^\circ \rangle = G$, so by simplicity, there is no nontrivial subgroup normalized by both $G_1^\circ$ and $G^\circ_y$. This then implies that $K$ does not contain $F^\circ(G_1)$ as otherwise it would imply that $F^\circ(G_1) = F^\circ(G_y)$.  

By Hrushovski, $K$ has corank $1$ or $2$ in $G^\circ_1$.  Combining this with our previous analysis, we find that the corank is equal to $2$ (hence $\rk K = 2$) and that $K$ contains a unique rank $1$ unipotent subgroup $W$ as well as some rank $1$ torus. In particular, $W$ is characteristic in $K$ so normal in $G^\circ_1$. We claim that  $W$ is also normal in $G^\circ_y$, and this will be our final contradiction.

Since $G^\circ_1/K$ has rank $2$, Hrushovski also implies that  $(G^\circ_{1})_y/K$ is a good torus, so  $(G^\circ_{1})_y$ contains a rank $2$ good torus $T$ of $G^\circ_y$. Since $W$ is characteristic in $K$, $T$ normalizes $W$. By the structure of $G^\circ_y$, $W\le F(G^\circ_y)$, implying that $F(G^\circ_y)$ also normalizes $W$. But then, $\langle T,F(G^\circ_y)  \rangle = G^\circ_y$ normalizes $W$ (as does $G^\circ_1$, hence $G$), a contradiction.
\end{proofclaim}

\begin{claim}\label{Claim.Prop.SimpleRank6.2PtStabOnG/D}
$G_1 = F^\circ(G_1) \rtimes G_{12}$. Thus, $(G,Y)$ is a split $2$-transitive group.
\end{claim}
\begin{proofclaim}
We begin by showing that $G_{1,2}^\circ$ is a good torus; assume not. 

As in the proof of Claim~\ref{Claim.Prop.SimpleRank6.2TransOnG/D}, no nontrivial subgroup is normalized by both $G_1^\circ$ and $G^\circ_2$, so $G_{1,2}^\circ\neq F^\circ(G_1^\circ)$. As we are assuming $G_{1,2}^\circ$ is not a good torus, we find that $G_{1,2}^\circ = WS$ for some rank $1$ unipotent group $W$ and some rank $1$ good torus $S$. Also, $W$ is normalized by $S$ since $G_{1,2}^\circ$ is solvable.

Consider $N:=N_G^\circ(W)$. It contains $\langle G_{1,2}^\circ, F^\circ(G_1^\circ),F^\circ(G_2^\circ)\rangle$, so $\rk N$ is at least $4$, hence equal to $4$. If $N$ is nonsolvable, the analysis of Claim~\ref{Claim.Prop.SimpleRank6.StructureM} could be applied to $N$ to force $W$ to be a torus, a contradiction. So $N$ is solvable, and by Claim~\ref{Claim.Prop.SimpleRank6.3Trans}, $N$ is nonnilpotent. Thus $F:= F^\circ(N)$ has rank $2$ or $3$.

First suppose $\rk F = 2$. Then $F^\circ(G_1^\circ)$ acts on the factors of $1\trianglelefteq W \trianglelefteq F$, which each have rank $1$. Any nontrivial action can be linearized, but $F^\circ(G_1^\circ)$ has no quotient isomorphic to the multiplicative group of a field. Thus, $F^\circ(G_1^\circ)$ acts trivially on each factor of the series, implying that $F^\circ(G_1^\circ) \le F$. Similarly, $F^\circ(G_2^\circ) \le F$, contradicting that $\rk F = 2$.

Next suppose $\rk F = 3$. Since $S$ does not commute with $F^\circ(G_1^\circ)$, $S$ is not contained in $F$, so $F=\langle F^\circ(G_1^\circ),F^\circ(G_2^\circ)\rangle$. Then $F^\circ(G_1^\circ)$ is corank $1$ in the nilpotent group $F$, so $F$ normalizes $F^\circ(G_1^\circ)$. Of course $F^\circ(G_1^\circ)$ is normal in $G_1^\circ$, and by rank considerations, $G_1^\circ$ is the connected normalizer in $G$ of $F^\circ(G_1^\circ)$. Thus, $F \le G_1^\circ$. However, we know the structure of $F$ (unipotent of rank $3$) and of $G_1^\circ$ (rank $4$ containing a rank $2$ good torus), so we have a contradiction.

Thus, $G_{1,2}^\circ$ is a good torus. We may  apply \cite[Lemma~3.7]{WiJ16} to see that $G_1$ and $G_{12}$ are connected, and then Claim~\ref{Claim.Prop.SimpleRank6.StructureD} completes the proof.
\end{proofclaim}

\begin{claim}
Contradiction.
\end{claim}
\begin{proofclaim}
We appeal to the theory of split $2$-transitive permutation groups (via the theory of Moufang sets). Since both factors in the splitting $G_1 = F^\circ(G_1) \rtimes G_{12}$ are abelian, we 
find that $G\cong \PSL_2(F)$  (see for example \cite[Corollary~1.2]{WiJ10}). However, proper definable subgroups of  $\PSL_2(F)$ are solvable while  $M$ is not.
\end{proofclaim}
This completes the proof of Proposition~\ref{Prop.SimpleRank6}.
\end{proof}\setcounter{claim}{0}

% % % % % % % % % % % % % % % % % % % %
% % % % % % % % % % % % % % % % % % % %
% SECTION 
% % % % % % % % % % % % % % % % % % % %
% % % % % % % % % % % % % % % % % % % %
\section{Proof of the main result}\label{S.Proof}

Throughout this final section, we adopt the hypotheses of the \hyperlink{h:theorem}{Theorem} as well as the additional notation presented in Section~\ref{s.Strategy}. 
In particular, $(G,X)$ is a generically sharply $t$-transitive permutation group of finite Morley rank, $(1,\ldots,t)\in X^t$ is a generic $t$-tuple, and  $r=\rk X$. By existing work of Hrushovski (see \cite[Theorem~11.98]{BoNe94}) and the present authors~\cite{AlWi18}, we may assume $r\ge 3$. 

We begin by addressing connectedness of point-stabilizers, which we use in our analysis of $G_{[t-1]}$ and $G_{[t-2]}$.

\begin{lemma}\label{lem.G[t-1]Connected}
If $t\ge 2$, then for every $0\le m < t$, the pointwise stabilizer $G_{[m]}$ of $1,\ldots,m$ is connected of rank $(t-m)\cdot r$, where we take $G_{[0]}$ to be $G$.
\end{lemma}
\begin{proof}   
Since $t\geq 2$, by \cite[Lemma 1.8(3)]{BoCh08}, the set $X$ is connected; this is a key point.

Now, as the action of $G$ is generically sharply $t$-transitive, $G_{[m]}$ acts generically sharply $(t-m)$-transitively, so if $\mathcal{O}\subseteq X^{t-m}$ is the generic orbit of $G_{[m]}$, then $G_{[m]}$ is in definable bijection with $\mathcal{O}$ (making essential use of  sharpness). Since $\mathcal{O}$ is generic in $X^{t-m}$, $\mathcal{O}$ has rank $(t-m)\cdot r$ as well as degree $1$ by our initial observation, so same is true of  $G_{[m]}$. 
\end{proof}

As indicated in Section~\ref{s.Strategy}, Proposition~\ref{prop.SymAltActionOnLGroup} will reduce our work to consideration of two cases: (1) $G_{[t-1]}$ is solvable or (2) $r = 3$ with $G_{[t-1]} \cong \PGL_2(F)$. We treat these two cases first and then our eventual proof of the \hyperlink{h:theorem}{Theorem} will simply tie things together.
We begin with the solvable case---this falls quickly thanks to Proposition~\ref{prop.Gen2TransAbelianPS}.

\begin{proposition}\label{prop.G[t-1]SolvableBound}
If $G_{[t-1]}$ is solvable, then $t\le r+2$.
\end{proposition}
\begin{proof}
Assume $G_{[t-1]}$ is solvable and, towards a contradiction, that $t \ge r+3$. Then $\Sigma_t \cong \Sym(t-1) \ge \Sym(r+2)$, and we consider the faithful action of $\Sigma_t$ on the rank $r$ group $G_{[t-1]}$ (see Section~\ref{s.Strategy}).

By \cite[Lemma~2.7]{AlWi24}, $G_{[t-1]}$ is abelian, and then \cite{CDW23} takes over---via the First Geometrisation Lemma followed by the Recognition Lemma---to show that $G_{[t-1]}$ is an elementary abelian $p$-group for some $p$ dividing $t-1$. However, since $G_{[t-1]}$ is abelian,  Proposition~\ref{prop.Gen2TransAbelianPS} applies to the action of $G_{[t-2]}$ on its generic orbit, and we find that $G_{[t-1]}$ has a quotient isomorphic to $L^\times$ for some algebraically closed field $L$, a contradiction.
\end{proof}

We now address the exceptional case of $r = 3$ and $G_{[t-1]} \cong \PGL_2(F)$. The goal is to eliminate this possibility, which we phrase as follows. 

\begin{proposition}\label{prop.Rank3HMustBeSolvable}
If  $r=3$ and  $t\ge 5$, then $G_{[t-1]}$ is solvable.
\end{proposition}

Our proof of Proposition~\ref{prop.Rank3HMustBeSolvable} will take some preparation.

\begin{lemma}\label{lem.Rank3StructureH-new}
Assume $r=3$, $t\ge 5$, and $G_{[t-1]}$ is nonsolvable. Then the following hold: 
\begin{enumerate}
    \item\label{lem.Rank3StructureH-new.item.G[t-1]PSL} $G_{[t-1]}\cong \PSL_2(K)$, for some rank $1$ field $K$, with  $S_t$ acting by inner automorphisms;
    \item $G_{[t-1]}$ is a maximal connected definable subgroup of $G_{[t-2]}$;
    \item\label{lem.Rank3StructureH-new.item.G[t-2]Quasisimple} $G_{[t-2]}$ is quasisimple with finite center.
\end{enumerate}

\end{lemma}
\begin{proof}
Set $B= G_{[t-2]}$ and $H = G_{[t-1]}$; recall from Lemma \ref{lem.G[t-1]Connected} that they are connected and of ranks $6$ and $3$ respectively. 

Given the classification of groups of rank $3$ (see \cite{ChG79,FrO18}) and the faithful action $S_t \ge \Sym(4)$ on $H$, we have that $H\cong \PGL_2(K)$, for some algebraically closed field $K$, with  $S_t$ acting by inner automorphisms.

We next show that $B$ is quasisimple with finite center. Suppose $N$ is a proper connected normal subgroup of $B$ and of maximal rank with respect to those properties; we aim to show $N=1$. Since $H$ is simple, $H\cap N = 1$, so $\rk N \le 3$ (since $\rk B = 6$). Now, $B/N$ is nonsolvable (as it embeds $H$), so by maximality of $N$, $B/N$ is quasisimple with finite center. And since there are no quasisimple groups of rank $4$ or $5$ that have involutions \cite{WiJ16,DeWi16}, $\rk B/N = \rk N = 3$.

We now have $HN = B$ and $N\cap H = 1$. Thus, $N$ acts regularly on the generic orbit of $B$ on $X$, and the action of $H$ on this orbit is isomorphic to the action of $H$ on $N$ by conjugation. Thus, $H$ acts generically transitively on $N$ by conjugation. If $N$ is solvable, then $N$ is abelian.  But then, \cite{BoDe15} implies that the action of $H$ on $N$ is the adjoint action, and this cannot be generically transitive since $H$ preserves the determinant. Otherwise, if $N$ is nonsolvable, then $N$ is quasisimple of rank $3$ so of the form $\pSL$. However, no such group supports a generically transitive automorphism group. Thus, $\rk N \neq 3$, so $B$ is indeed quasisimple with finite center.

Finally, we show that $H$ is a maximal connected definable subgroup of $B$. Indeed, if $M$ is a  proper connected definable subgroup of $B$ that properly contains $H$, then $M$ is nonsolvable and $\rk M \ge 4$. Moving to the quotient of $B$ by its center, we have a contradiction to Proposition \ref{Prop.SimpleRank6}.
\end{proof}

Our contradictory assumption that $G_{[t-1]}$ is nonsolvable forces $G_{[t-2]}$ to be quasisimple with a finite center and  many involutions (due to the structure of $G_{[t-1]}$), and this creates serious tension with the Algebraicity Conjecture. Since $G_{[t-2]}$ has rank $6$, the conjecture implies---as seen in Lemma \ref{smallalgebraicgroups} below---that $G_{[t-2]}$ should be of the form $\PSL_2$ so should certainly \emph{not} have a subgroup like $G_{[t-1]}$. This is exactly our contradiction when $G_{[t-2]}$ is of even type thanks to \cite{ABC08}. In odd type, nonalgebraic $G_{[t-2]}$ remains a possibility, but we are able to sufficiently pin down its structure---this time thanks to Proposition~\ref{Prop.SimpleRank6} and \cite{DeJa16}---to again reach a contradiction. 

\begin{lemma}\label{smallalgebraicgroups}
The only infinite simple algebraic groups of Morley rank at most $6$ are $\PSL_2(K)$ over a field $K$ of Morley rank $1$ or $2$. 
\end{lemma}
\begin{proof}
Let $G$ be an infinite simple algebraic group over a field $K$ of finite Morley rank. By a theorem of Macintyre, $K$ is algebraically
closed, and thus, the structure theory of simple algebraic groups over algebraically closed fields applies. According to this, the dimension
of $G$ over the pure field structure is $2u+l$ where $u$ is the dimension of the unipotent subgroup generated by a fixed set of positive
roots and $l$ is the Lie rank of $G$ (see for example Section 28.5 of \cite{HuJ75}). Moreover, $u\geq l\geq 1$. 
It follows that $\rk(G)=(2u+l)\rk(K)\geq 3\rk(K)\geq 3$.
If $l>1$ then the smallest value of $u$ is $3$ (see for example Section 9.3 of \cite{HuJ72}). Thus, $l=1$ since $\rk(G)\leq6$.
In this case, the only possibility for $G$ is $\PSL_2(K)$ (Corollary 32.3 of \cite{HuJ75}), and $\rk(K)$ is either $1$ or $2$. 
\end{proof}

\begin{corollary}\label{cor.Boddtype}
Under the assumptions of Lemma~\ref{lem.Rank3StructureH-new}, $G_{[t-2]}$ is of odd type. 
\end{corollary}
\begin{proof}
Set $B=G_{[t-2]}$. Then $B$ is  quasisimple with finite center according to Lemma \ref{lem.Rank3StructureH-new}\eqref{lem.Rank3StructureH-new.item.G[t-2]Quasisimple}. If $B$ contains a nontrivial unipotent $2$-subgroup, so does $B/Z(B)$. By the classification of the simple groups of finite Morley rank of even and mixed type (\cite{ABC08}), $B/Z(B)$ is a simple linear algebraic group over an algebraically closed field of characteristic $2$. Lemma \ref{smallalgebraicgroups} then implies  $B/Z(B)\cong\PSL_2(K)$ with $K$ algebraically closed, but $G_{[t-1]}$ is not solvable by Lemma \ref{lem.Rank3StructureH-new}\eqref{lem.Rank3StructureH-new.item.G[t-1]PSL}, violating the structure of $\PSL_2(K)$. 
\end{proof}

The tension in odd type remains. Things will begin to give as we invoke Proposition~\ref{Prop.SimpleRank6} and ultimately Fact~\ref{fact.DJ16} to show that if $G_{[t-2]}$ is not of the form $\PSL_2$ then it is still---in some sense---close (being of type \cibo$_2$). In any case, the structure of $G_{[t-2]}$ becomes rather transparent and, with a little more effort, quickly falls apart.

The next lemma provides the bridge from Proposition~\ref{Prop.SimpleRank6} to Fact~\ref{fact.DJ16}; the proof of Proposition~\ref{prop.Rank3HMustBeSolvable} will follow.

\begin{lemma}\label{lem.Rank3HCentralizerInvolution}
Under the assumptions of Lemma~\ref{lem.Rank3StructureH-new}, $C^\circ_{G_{[t-2]}}(i)$ is solvable for each involution $i\in G_{[t-2]}$.
\end{lemma}
\begin{proof}
Let $B= G_{[t-2]}$ and $C = C^\circ_{B}(i)$. By \cite[Proposition 5.15]{BoCh08}, $i\in C$.
Suppose now that $C$ is not solvable. By Proposition~\ref{Prop.SimpleRank6}, $\rk C \le 3$; this is a crucial point that eliminates the possibility of $C$ containing a copy of $\PSL_2$. 

The classification of nonsolvable groups of rank at most $3$ (and  that $i\in C$) implies that $C\cong \SL_2(F)$ with $F$ of characteristic not $2$. But then the Sylow $2$-subgroup of $B$ matches that of $\SL_2(F)$, contrary to the fact that the Sylow $2$-subgroup of $G_{[t-1]}\le B$ matches that of $\PSL_2(F)$ by Lemma~\ref{lem.Rank3StructureH-new}\eqref{lem.Rank3StructureH-new.item.G[t-1]PSL}.
\end{proof}

\begin{proof}[Proof of Proposition~\ref{prop.Rank3HMustBeSolvable}]
We  work with a counterexample; assume:  $r=3$, $t\ge 5$, and $H = G_{[t-1]}$ is nonsolvable. Set $B= G_{[t-2]}$. By  Lemma~\ref{lem.Rank3StructureH-new} and Corollary~\ref{cor.Boddtype} together with Proposition~\ref{Prop.SimpleRank6}, $B$ is an $N_\circ^\circ$-group of odd type, and by Lemma~\ref{lem.Rank3HCentralizerInvolution}, we may apply Fact~\ref{fact.DJ16}. The presence and structure of $H$ given in Lemma~\ref{lem.Rank3StructureH-new} implies that we are in type \cibo$_2$. 

Consider  $C = C_B^\circ(i)$ for any $i\in I(B)$.  Then $C$ is an abelian Borel subgroup of $B$ and $\rk C = 2$. We claim that $C$ is maximal connected in $B$. Indeed, if not, then $C < M$ with $M$ connected of rank $3$ or $4$, and as $C$ is a \emph{Borel} subgroup of $B$, $M$ is nonsolvable. By Proposition~\ref{Prop.SimpleRank6}, $M$ must have rank $3$, but no  nonsolvable rank $3$ group has a rank $2$ abelian subgroup. Thus $C$ is indeed maximal connected in $B$. 

Next, notice that if $i \in H$, then $T:= (C\cap H)^\circ$ is a rank one good torus, which is  the definable closure of a $2$-torus in $H$. In particular, $T$ is a characteristic rank $1$ subgroup of $C$. And as involutions are conjugate in $B$ (as a consequence of having type \cibo$_2$), $C_B^\circ(i)$ has such a characteristic rank $1$ subgroup, for every choice of $i\in I(B)$.

We  apply our observations to $C = C_B^\circ(\omega)$ for $\omega = (t-1, t)$. Note that $H^\omega \cap H = 1$, so $H \cap C = 1$. We  study the  action of $\Sigma_{t-1,t} \cong\Sym(t-2)\ge \Sym(3)$ on $C$ and on $H$ to show, for a $3$-cycle $\gamma \in \Sigma_{t-1,t}$, that  $C_B^\circ(\gamma)$ contains $C$ and has rank at least $3$, contradicting that  $C$ is maximal connected in $B$. 

Let $T$ be the definable closure of a $2$-torus in $C$. Then $\Sigma_{t-1,t}$ acts on the rank $1$ factors of the series $1\le T\le C$. Thus, $[\gamma,C] \le T$, and $[\gamma,\gamma,C] = 1$. If $[\gamma,C] \ne 1$, then (by usual quadratic arguments) $[\gamma,C]$ is an elementary abelian $3$-group, which also must be equal to $T$, a contradiction. We conclude that $[\gamma,C] = 1$, so $C\le C_B^\circ(\gamma)$.  

By Lemma~\ref{lem.Rank3StructureH-new}\eqref{lem.Rank3StructureH-new.item.G[t-1]PSL},   $\gamma$ acts on $H$ as an inner automorphism, so $\gamma$  centralizes a rank $1$ connected subgroup $A$ of $H$. As observed above, $H \cap C = 1$, so $\langle C, A\rangle$ is a connected subgroup of $C_B^\circ(\gamma)$ of rank at least $3$. Additionally,  $\Sigma_{t-1,t}$ acts faithfully on $B$ (see Section~\ref{s.Strategy}), so $C_B^\circ(\gamma)$ is  proper subgroup of $B$ that properly contains $C$, a final contradiction. 
\end{proof}

Now, we complete the proof of the main theorem of the paper.

\begin{proof}[Proof of the \hyperlink{h:theorem}{Theorem}]
We work towards a contradiction under the assumption that $t\geq r+3$. As we noted at the beginning of this section, we may assume that $r\geq3$. Let $H=G_{[t-1]}$. By Lemma \ref{lem.G[t-1]Connected}, $H$ is connected of Morley rank $r$.

By Proposition~\ref{prop.G[t-1]SolvableBound}, $H$ is not solvable, and then invoking Proposition~\ref{prop.Rank3HMustBeSolvable}, we have that $r\ge 4$. %Recall from Section~\ref{s.Strategy} that 
We have seen that $H$ is connected of rank $r$ and $\Sigma_t\cong\Sym(t-1)$ acts faithfully on $H$. As $t-1\ge r+2 \ge 6$, Proposition~\ref{prop.SymAltActionOnLGroup} implies that 
$t-1 \le \rk(H) = r$, a contradiction.
\end{proof}

% % % % % % % % % % % % % % % % % % % %
% % % % % % % % % % % % % % % % % % % %
% SECTION 
% % % % % % % % % % % % % % % % % % % %
% % % % % % % % % % % % % % % % % % % %
\section*{Acknowledgements}
We warmly thank the anonymous referee for a  careful reading of the article and numerous helpful suggestions. Among other things, their remarks led us to clarify a number of points and  to identify a small case that was overlooked in an earlier version of Proposition~\ref{prop.SymAltActionOnLGroup}.

The work of the second author was partially supported by the National Science Foundation under grant No.~DMS-1954127.
\bibliographystyle{alpha}
\bibliography{GenericTransitivity}
\end{document}